\documentclass[12pt,english]{amsart}
\usepackage[T1]{fontenc}
\usepackage[latin9]{inputenc}
\usepackage{geometry}
\usepackage{babel}
\usepackage{array}
\usepackage{amsthm}
\usepackage{amsbsy}
\usepackage{amstext}
\usepackage{amssymb}
\usepackage{esint}
\usepackage[unicode=true,pdfusetitle,
 bookmarks=true,bookmarksnumbered=false,bookmarksopen=false,
 breaklinks=false,pdfborder={0 0 1},backref=false,colorlinks=false]
 {hyperref}

\makeatletter

\providecommand{\tabularnewline}{\\}

 \theoremstyle{definition}
 \newtheorem*{defn*}{\protect\definitionname}
\theoremstyle{plain}
\newtheorem{thm}{\protect\theoremname}

\usepackage{subfig}

\usepackage{tikz}
\usetikzlibrary{decorations.markings, arrows}

\newlength{\myArraycolsep}
\@ifundefined{showcaptionsetup}{}{%
 \PassOptionsToPackage{caption=false}{subfig}}
\usepackage{subfig}
\makeatother

  \providecommand{\definitionname}{Definition}
\providecommand{\theoremname}{Theorem}

\begin{document}

\title[Isospectrality via eigenderivative transplantation]{Changing gears:\\
Isospectrality via eigenderivative transplantation}

\author{Peter Doyle and Peter Herbrich}
\begin{abstract}
We introduce a new method for constructing isospectral quantum graphs
that is based on transplanting derivatives of eigenfunctions. We also
present simple digraphs with the same reversing zeta function, which
generalizes the Bartholdi zeta function to digraphs.
\end{abstract}

\subjclass[2000]{34B45, 05C50, 58J53}

\keywords{Isospectrality, Quantum graphs, Markov chains, Digraphs, Zeta functions}

\address{Department of Mathematics, Dartmouth College, Hanover, NH, USA}

\email{peter.g.doyle@dartmouth.edu, peter.herbrich@dartmouth.edu}

\thanks{We thank Ram Band who conjectured that the graphs in Figure~\ref{fig:Balanced_three_gears}
feature a new kind of transplantation.}

\maketitle
\begin{minipage}[t]{1\columnwidth}%
\global\long\def\graph{G}

\global\long\def\other#1{\tilde{#1}}

\global\long\def\nrVertices#1{n_{#1}}
\global\long\def\vertices{V}

\global\long\def\vertex#1{v_{#1}}

\global\long\def\vertexOther#1{\other v_{#1}}

\global\long\def\nrEdges#1{m_{#1}}
\global\long\def\edges{E}
\global\long\def\edge#1{e_{#1}}

\global\long\def\side#1{p_{#1}}

\global\long\def\sideOther#1{\other p_{#1}}

\global\long\def\tooth#1{t_{#1}}

\global\long\def\toothOther#1{\other t_{#1}}

\global\long\def\lengths{\boldsymbol{l}}

\global\long\def\length#1{l_{#1}}

\global\long\def\lengtha{a}

\global\long\def\lengthb{b}

\global\long\def\lengthc{c}

\global\long\def\lengthd{d}

\global\long\def\weight{w}

\global\long\def\eigenvalue{\lambda}

\global\long\def\eigenvalueMarkov{\mu}

\global\long\def\identityMatrix#1{I_{#1}}

\global\long\def\allOnesMatrix#1{J_{#1}}

\global\long\def\allOnesMatrixOther#1{\other J_{#1}}

\global\long\def\degreeMatrix#1{D_{#1}}

\global\long\def\adjacencyMatrix#1{A_{#1}}

\global\long\def\indegreeMatrix#1{D_{#1}^{\mathrm{in}}}

\global\long\def\outdegreeMatrix#1{D_{#1}^{\mathrm{out}}}

\global\long\def\varDout{\gamma}

\global\long\def\varDin{\delta}

\global\long\def\varA{\alpha}

\global\long\def\varAT{\beta}

\global\long\def\generalizedCharacteristicPolynomial#1{\eta_{#1}}

\global\long\def\Trace{\mathrm{Tr}}

\global\long\def\Determinant{\mathrm{Det}}

\global\long\def\CharacteristicPolynomial#1{\chi_{#1}}
\end{minipage}

\thispagestyle{empty}

\section{Introduction}

Quantum graphs are singular one-dimensional manifolds equipped with
self-adjoint operators, whose spectra allow for explicit computations.
Following~\cite{BolteEndres2009}, let $\graph=(\vertices,\edges,\lengths)$
be a finite metric graph with vertices $\vertices=\{\vertex 1,\vertex 2,\ldots,\vertex{\nrVertices{}}\}$,
edges $\edges=\{\edge 1,\edge 2,\ldots,\edge{\nrEdges{}}\}$, and
edge lengths $\lengths=(\length 1,\length 2,\ldots,\length{\nrEdges{}})\in\mathbb{R}_{+}^{\nrEdges{}}$.
In particular, $\edge i$ is parameterized by $x_{i}\in[0,\length i]$,
which determines an orientation of the edges. However, the differential
operator of interest acts as $\Delta_{i}=-\frac{\partial^{2}}{\partial x_{i}^{2}}$
on $C_{0}^{\infty}(0,\length i)$ and is thus invariant under transformations
of the form $x_{i}'=\length i-x_{i}$, which allows to regard edges
as undirected. The self-adjoint extensions of the corresponding symmetric
operator with initial domain $\oplus_{i=1}^{\nrEdges{}}H_{0}^{2}(0,\length i)$
can be parameterized in terms of boundary conditions at the vertices~\cite{KostrykinSchrader1999},
which turn the metric graph into a quantum graph. Unless otherwise
stated, we consider Kirchhoff-Neumann conditions, which require that
for each vertex~$\vertex i$, functions on adjacent edges take the
same value at~$\vertex i$, and the sum of their outgoing derivatives
at~$\vertex i$ vanishes. In particular, leaf vertices carry Neumann
boundary conditions. It is well-known that finite quantum graphs have
discrete spectrum~\cite{BerkolaikoKuchment2013}. Quantum graphs
are called isospectral if their spectra coincide, including multiplicities.

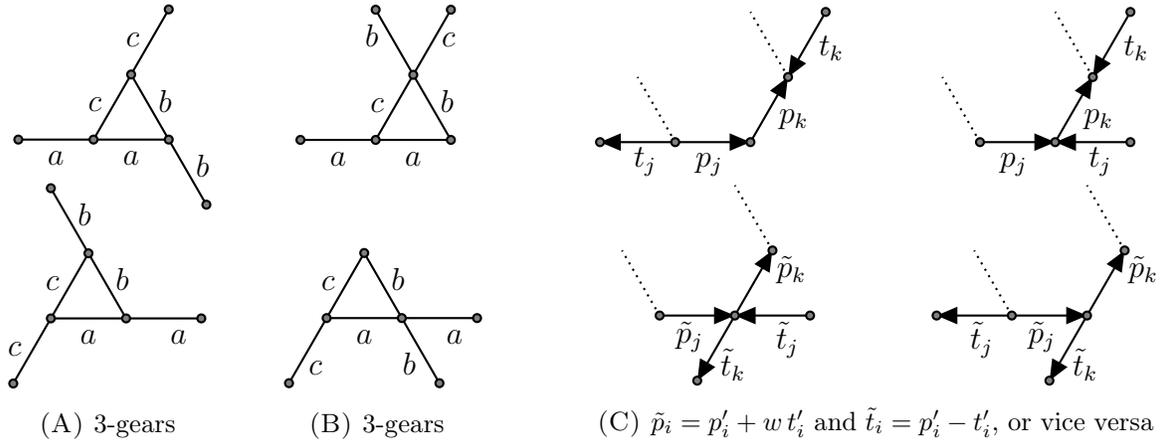
\begin{figure}
\begin{centering}
\subfloat[$3$-gears\label{fig:Balanced_three_gears}]{\protect\begin{centering}
\begin{minipage}[b][30mm]{0.23\columnwidth}%
\protect\begin{center}
\begin{tabular}{c}
\vspace{-3mm}
\begin{tikzpicture}[>=triangle 45,thick, decoration =
       {markings,mark=at position 1 with {\arrow{>}}}]
 \tikzstyle{every node} = [circle, draw, fill=black!50,
                           inner sep=1pt, minimum width=3pt]
 \tikzstyle{txt} = [rectangle, draw=none, fill=white, inner sep=-1pt]

 \draw ( 0, 0) --   (-1, 0);
 \draw ( 0, 0) --   ( 1, 0);
 \draw ( 1, 0) --   ( 60: 1);
 \draw ( 1, 0) -- ++(-60: 1);
 \draw (60: 1) --   (0, 0);
 \draw (60: 1) --   (60:2);

 \draw (-1,0) node{} (0,0) node{} (1,0) node{} (1,0)++(-60: 1) node{} (60:1) node{} (60:2) node{};

 \draw ( 0.5,-0.25) node[txt] {$\lengtha$};
 \draw (-0.5,-0.25) node[txt] {$\lengtha$};
 \draw ( 1.25,0)++(120:0.6) node[txt] {$\lengthb$};
 \draw ( 1.25,0)++(-60:0.4) node[txt] {$\lengthb$};
 \draw (-0.25,0)++( 60:0.54) node[txt] {$\lengthc$};
 \draw (-0.25,0)++( 60:1.54) node[txt] {$\lengthc$};

\end{tikzpicture}\vspace{-2mm}
\tabularnewline
\begin{tikzpicture}[>=triangle 45,thick, decoration =
       {markings,mark=at position 1 with {\arrow{>}}}]
 \tikzstyle{every node} = [circle, draw, fill=black!50,
                           inner sep=1pt, minimum width=3pt]
 \tikzstyle{txt} = [rectangle, draw=none, fill=white, inner sep=-1pt]

 \draw ( 0,0) --   (   1,0);
 \draw ( 1,0) --   (   2,0);
 \draw ( 1,0) --   (  60:1);
 \draw (60:1) -- ++( 120:1);
 \draw (60:1) --   (  0, 0);
 \draw ( 0,0) --   (-120:1);

 \draw (0,0) node{} (1,0) node{} (60:1) node{}
  (2,0) node{} (1,0)++(120:2) node{} (-120:1) node{};

 \draw ( 0.5,-0.25) node[txt] {$\lengtha$};
 \draw ( 1.7,-0.25) node[txt] {$\lengtha$};
 \draw ( 1.25,0)++( 120:0.6 ) node[txt] {$\lengthb$};
 \draw ( 1.25,0)++( 120:1.6 ) node[txt] {$\lengthb$};
 \draw (-0.25,0)++(  60:0.54) node[txt] {$\lengthc$};
 \draw (-0.25,0)++(-120:0.46) node[txt] {$\lengthc$};

\end{tikzpicture}\tabularnewline
\end{tabular}\protect
\par\end{center}%
\end{minipage}\protect
\par\end{centering}

}\subfloat[$3$-gears\label{fig:Unbalanced_three_gears}]{\protect\begin{centering}
\begin{minipage}[b][30mm]{0.17\columnwidth}%
\protect\begin{center}
\begin{tabular}{c}
\begin{tikzpicture}[>=triangle 45,thick, decoration =
       {markings,mark=at position 1 with {\arrow{>}}}]
 \tikzstyle{every node} = [circle, draw, fill=black!50,
                           inner sep=1pt, minimum width=3pt]
 \tikzstyle{txt} = [rectangle, draw=none, fill=white, inner sep=-1pt]

 \draw ( 0, 0) --   (-1, 0);
 \draw ( 0, 0) --   ( 1, 0);
 \draw ( 1, 0) --   ( 60: 1);
 \draw (60: 1) -- ++(120: 1);
 \draw (60: 1) --   (0, 0);
 \draw (60: 1) --   (60:2);

 \draw (-1,0) node{} (0,0) node{} (1,0) node{} (60:1)++(120: 1) node{} (60:1) node{} (60:2) node{};

 \draw ( 0.5,-0.25) node[txt] {$\lengtha$};
 \draw (-0.5,-0.25) node[txt] {$\lengtha$};
 \draw ( 0.75,0)++(120:1.6) node[txt] {$\lengthb$};
 \draw ( 1.25,0)++(120:0.6) node[txt] {$\lengthb$};
 \draw (-0.25,0)++( 60:0.54) node[txt] {$\lengthc$};
 \draw ( 0.2,0)++( 60:1.54) node[txt] {$\lengthc$};

\end{tikzpicture}\tabularnewline
\tabularnewline
\tabularnewline
\begin{tikzpicture}[>=triangle 45,thick, decoration =
       {markings,mark=at position 1 with {\arrow{>}}}]
 \tikzstyle{every node} = [circle, draw, fill=black!50,
                           inner sep=1pt, minimum width=3pt]
 \tikzstyle{txt} = [rectangle, draw=none, fill=white, inner sep=-1pt]

 \draw ( 0,0) --   (   1,0);
 \draw ( 1,0) --   (   2,0);
 \draw ( 1,0) --   (  60:1);
 \draw ( 1,0) -- ++( -60:1);
 \draw (60:1) --   (  0, 0);
 \draw ( 0,0) --   (-120:1);

 \draw (0,0) node{} (1,0) node{} (60:1) node{}
       (2,0) node{} (1,0)++(-60:1) node{} (-120:1) node{};

 \draw ( 0.5,-0.25) node[txt] {$\lengtha$};
 \draw ( 1.7,-0.25) node[txt] {$\lengtha$};
 \draw ( 1.25,0)++(120:0.6) node[txt] {$\lengthb$};
 \draw ( 1.1,-0.6)            node[txt] {$\lengthb$};
 \draw (-0.25,0)++(  60:0.54) node[txt] {$\lengthc$};
 \draw (-0.15,-0.66) node[txt] {$\lengthc$};

\end{tikzpicture}\tabularnewline
\end{tabular}\protect
\par\end{center}%
\end{minipage}\protect
\par\end{centering}

}\subfloat[$\protect\sideOther i=\protect\side i'+\protect\weight\,\protect\tooth i'$
and $\protect\toothOther i=\protect\side i'-\protect\tooth i'$, or
vice versa\label{fig:Polygon_vertices}]{\protect\begin{centering}
\begin{minipage}[b][30mm]{0.6\columnwidth}%
\protect\begin{center}
\begin{tabular}{c>{\centering}p{5mm}c}
\begin{tikzpicture}[>=triangle 45,thick, decoration =
       {markings,mark=at position 1 with {\arrow{>}}}]
 \tikzstyle{every node} = [circle, draw, fill=black!50,
                           inner sep=1pt, minimum width=3pt]
 \tikzstyle{txt} = [rectangle, draw=none, fill=white, inner sep=-1pt]

 \draw[postaction={decorate}] (-1, 0) -- (-2, 0);
 \draw[postaction={decorate}] (-1, 0) -- ( 0, 0);
 \draw[postaction={decorate}] (60: 0) -- (60: 1);
 \draw[postaction={decorate}] (60: 2) -- (60: 1);
 \draw[dotted] (-1, 0) -- ++(120:1);
 \draw[dotted] (60: 1) -- ++(120:1);

 \draw (-2,0) node{} (-1,0) node{} (0,0) node{} (60:1) node{} (60:2) node{};

 \draw (-0.55,-0.3) node[txt] {$\side j$};
 \draw (-1.4,-0.26) node[txt] {$\tooth j$};

 \draw (60: 0.35)++(0.4,0) node[txt] {$\side k$};
 \draw (60: 1.45)++(0.35,0) node[txt] {$\tooth k$};
\end{tikzpicture} &  & \begin{tikzpicture}[>=triangle 45,thick, decoration =
       {markings,mark=at position 1 with {\arrow{>}}}]
 \tikzstyle{every node} = [circle, draw, fill=black!50,
                           inner sep=1pt, minimum width=3pt]
 \tikzstyle{txt} = [rectangle, draw=none, fill=white, inner sep=-1pt]

 \draw[postaction={decorate}] (-1, 0) -- ( 0, 0);
 \draw[postaction={decorate}] ( 1, 0) -- ( 0, 0);
 \draw[postaction={decorate}] (60: 0) -- (60: 1);
 \draw[postaction={decorate}] (60: 2) -- (60: 1);
 \draw[dotted] (-1, 0) -- ++(120:1);
 \draw[dotted] (60: 1) -- ++(120:1);

 \draw (-1,0) node{} (0,0) node{} (1,0) node{} (60:1) node{} (60:2) node{};

 \draw (-0.55,-0.3) node[txt] {$\side j$};
 \draw ( 0.6,-0.26) node[txt] {$\tooth j$};

 \draw (60: 0.35)++(0.4,0) node[txt] {$\side k$};
 \draw (60: 1.45)++(0.35,0) node[txt] {$\tooth k$};
\end{tikzpicture}\tabularnewline
\begin{tikzpicture}[>=triangle 45,thick, decoration =
       {markings,mark=at position 1 with {\arrow{>}}}]
 \tikzstyle{every node} = [circle, draw, fill=black!50,
                           inner sep=1pt, minimum width=3pt]
 \tikzstyle{txt} = [rectangle, draw=none, fill=white, inner sep=-1pt]

 \draw[postaction={decorate}] (-1, 0) -- ( 0, 0);
 \draw[postaction={decorate}] ( 1, 0) -- ( 0, 0);
 \draw[postaction={decorate}] (60: 0) -- (60:-1);
 \draw[postaction={decorate}] (60: 0) -- (60: 1);
 \draw[dotted] (-1, 0) -- ++(120:1);
 \draw[dotted] (60: 1) -- ++(120:1);

 \draw (-1,0) node{} (0,0) node{} (1,0) node{} (60:-1) node{} (60:1) node{};

 \draw (-0.6,-0.28) node[txt] {$\sideOther j$};
 \draw ( 0.7,-0.32) node[txt] {$\toothOther j$};

 \draw (60: 0.7)++(0.4,0) node[txt] {$\sideOther k$};
 \draw (60:-0.75)++(0.35,0) node[txt] {$\toothOther k$};

\end{tikzpicture} &  & \begin{tikzpicture}[>=triangle 45,thick, decoration =
       {markings,mark=at position 1 with {\arrow{>}}}]
 \tikzstyle{every node} = [circle, draw, fill=black!50,
                           inner sep=1pt, minimum width=3pt]
 \tikzstyle{txt} = [rectangle, draw=none, fill=white, inner sep=-1pt]

 \draw[postaction={decorate}] (-1, 0) -- (-2, 0);
 \draw[postaction={decorate}] (-1, 0) -- ( 0, 0);
 \draw[postaction={decorate}] (60: 0) -- (60:-1);
 \draw[postaction={decorate}] (60: 0) -- (60: 1);
 \draw[dotted] (-1, 0) -- ++(120:1);
 \draw[dotted] (60: 1) -- ++(120:1);

 \draw (-2,0) node{} (-1,0) node{} (0,0) node{} (60:-1) node{} (60:1) node{};

 \draw (-0.6,-0.28) node[txt] {$\sideOther j$};
 \draw (-1.4,-0.32) node[txt] {$\toothOther j$};

 \draw (60: 0.7)++(0.4,0) node[txt] {$\sideOther k$};
 \draw (60:-0.75)++(0.35,0) node[txt] {$\toothOther k$};

\end{tikzpicture}\tabularnewline
\end{tabular}\protect
\par\end{center}%
\end{minipage}\protect
\par\end{centering}

}
\par\end{centering}

\protect\caption{Eigenderivative transplantation. At the leaf vertices, we impose Neumann
conditions. At the polygon vertices, we require continuity as well
as $\protect\side j'(\protect\length j)+\protect\weight\,\protect\tooth j'(\protect\length j)=\protect\side k'(0)+\protect\weight\,\protect\tooth k'(0)$
and $\protect\sideOther j'(\protect\length j)+\protect\weight\,\protect\toothOther j'(\protect\length j)=\protect\sideOther k'(0)+\protect\weight\,\protect\toothOther k'(0).$\label{fig:Eigenderivative_transplantation}}
\end{figure}

Recently, it was discovered~\cite{OrenBand2012} that the graphs
in Figure~\ref{fig:Balanced_three_gears} have the same set of eigenvalues
when viewed as either edge-weighted combinatorial graphs or quantum
graphs with edge weights or lengths $(a,a,b,b,c,c)$, respectively.
These graphs first appeared in~\cite{McDonaldMeyers2003} where a
certain line graph construction was applied to the $7_{1}$ pair of
Dirichlet isospectral planar domains in~\cite{BuserConwayDoyleSemmler1994}.
The second author~\cite{Herbrich2015} has generalized this construction
to manifolds with mixed Dirichlet-Neumann boundary conditions and
revealed its connection to the graph-theoretic characterization of
the famous Sunada method~\cite{Sunada1985} given in~\cite{Herbrich2011}.

It is worth mentioning that the main arguments in~\cite{OrenBand2012}
are based on the widespread misconception that eigenfunctions on quantum
graphs with Kirchhoff-Neumann conditions are determined by their values
at the vertices~\cite{KottosSmilansky1999,ShapiraSmilansky2006,BandSmilansky2007}.
For example, if $\lengtha$, $\lengthb$, and $\lengthc$ are integer
multiples of some $r>0$, then each of the graphs in Figures~\ref{fig:Balanced_three_gears}
and~\ref{fig:Unbalanced_three_gears} has countably many eigenfunctions
which are supported on its central cycle of length $a+b+c$ and vanish
at all vertices. Instead, \cite{OrenBand2012}~asserts that if $f_{i}\colon[0,\length i]\to\mathbb{R}$
denotes the restriction of an eigenfunction to the edge $\edge i$
with eigenvalue $\eigenvalue\geq0$, then 
\begin{equation}
f_{i}(x_{i})=\frac{1}{\sin(\sqrt{\eigenvalue}\length i)}\left(f_{i}(0)\sin(\sqrt{\eigenvalue}(\length i-x_{i}))+f_{i}(\length i)\sin(\sqrt{\eigenvalue}x_{i})\right).\label{eq:Eigenfunction_in_terms_of_its_boundary_values}
\end{equation}
Similarly, \cite{KottosSmilansky1999} suggests to ignore edges $\edge i$
satisfying $\sqrt{\eigenvalue}\length i=k\pi$ for some $k\in\mathbb{Z}_{+}$
when determining whether $\eigenvalue$ is an eigenvalue, while such
edges contribute the boundary conditions
\begin{equation}
f_{i}(\length i)=(-1)^{k}f_{i}(0)\qquad\text{and}\qquad f_{i}'(\length i)=(-1)^{k}f_{i}'(0).\label{eq:Boundary_condition_for_Dirichlet_eigenvalues}
\end{equation}

In Section~\ref{sec:Eigenderivative_Transplantation}, we introduce
a method for constructing isospectral quantum graphs that avoids explicit
computations and thus bypasses the above-mentioned shortcomings. It
is a derivative of Buser's transplantation method~\cite{Buser1986},
which itself can be viewed as the combinatorial incarnation of the
Sunada method~\cite{Sunada1985}. In contrast to the latter, our
method can produce pairs of quantum graphs without common covers.
In Section~\ref{sec:Combinatorial_eigenderivative_transplantation},
we relate our method to random walks on combinatorial graphs. In addition,
we present related pairs of non-regular simple digraphs that have
the same reversing zeta function as introduced in~\cite{Herbrich2014}.
Thus, these digraphs exhibit a noteworthy degree of spectral indistinguishability.

\section{Gear graphs and eigenderivative transplantation \label{sec:Eigenderivative_Transplantation}}

\global\long\def\nrTeeth{n}

\begin{figure}
\begin{centering}
\subfloat[Derived from $3$-gears with $\protect\weight=3/2.$\label{fig:Derived_graphs_weight_three_halfs}]{\protect\centering{}%
\begin{minipage}[b][27mm]{0.48\columnwidth}%
\protect\begin{center}
 \begin{tikzpicture}[thick]
 \tikzstyle{every node} = [circle, draw, fill=black!50, 
                           inner sep=1pt, minimum width=3pt]
 \tikzstyle{txt} = [rectangle, draw=none, fill=white, inner sep=-1pt]

 \draw ( 0:0) to (170:1) node{};
 \draw ( 0:0) to (180:1) node{};
 \draw ( 0:0) to (190:1) node{};

 \draw (60:1) to             ++(120:1) node{};
 \draw (60:1) to[bend left]  ++(120:1) node{};
 \draw (60:1) to[bend right] ++(120:1) node{};

 \draw (60:1) to             (60:2) node{};
 \draw (60:1) to[bend left]  (60:2) node{};
 \draw (60:1) to[bend right] (60:2) node{};

 \draw (0:0) to[bend left] (0:1);
 \draw (0:0) to[bend right] (0:1);

 \draw (0:1) to[bend left] (60:1);
 \draw (0:1) to[bend right] (60:1);

 \draw (60:1) to[bend left] (0:0);
 \draw (60:1) to[bend right] (0:0);

 \draw (0:0) node{} (0:1) node{} (60:1) node{};

 \draw ( 0:0.5)          node[txt]{$\lengtha$};
 \draw (60:0.5)++(0.52,0) node[txt]{$\lengthb$};
 \draw (60:0.5)          node[txt]{$\lengthc$};

 \draw (-0.5,-0.3) node[txt]{$\lengtha$};
 \draw (-0.14, 1.3) node[txt]{$\lengthb$};
 \draw ( 1.1, 1.25) node[txt]{$\lengthc$};

\end{tikzpicture}\hspace*{5mm} \begin{tikzpicture}[thick]
 \tikzstyle{every node} = [circle, draw, fill=black!50, 
                           inner sep=1pt, minimum width=3pt]
 \tikzstyle{txt} = [rectangle, draw=none, fill=white, inner sep=-1pt]

 \draw ( 0:1) to ++(-10:1) node{};
 \draw ( 0:1) to ++(  0:1) node{};
 \draw ( 0:1) to ++( 10:1) node{};

 \draw ( 0:1) to             ++(-60:1) node{};
 \draw ( 0:1) to[bend left]  ++(-60:1) node{};
 \draw ( 0:1) to[bend right] ++(-60:1) node{};

 \draw (0:0) to             (-120:1) node{};
 \draw (0:0) to[bend left]  (-120:1) node{};
 \draw (0:0) to[bend right] (-120:1) node{};

 \draw (0:0) to[bend left] (0:1);
 \draw (0:0) to[bend right] (0:1);

 \draw (0:1) to[bend left] (60:1);
 \draw (0:1) to[bend right] (60:1);

 \draw (60:1) to[bend left] (0:0);
 \draw (60:1) to[bend right] (0:0);

 \draw (0:0) node{} (0:1) node{} (60:1) node{};

 \draw ( 0:0.5)          node[txt]{$\lengtha$};
 \draw (60:0.5)++(0.5,0) node[txt]{$\lengthb$};
 \draw (60:0.5)          node[txt]{$\lengthc$};

 \draw (1.7,-0.3) node[txt]{$\lengtha$};
 \draw (0.95,-0.6) node[txt]{$\lengthb$};
 \draw (0.0,-0.65) node[txt]{$\lengthc$};

\end{tikzpicture}\protect
\par\end{center}%
\end{minipage}\protect}\subfloat[Derived from $4$-gears with ${\displaystyle \protect\weight=2.}$]{\protect\begin{centering}
\begin{minipage}[b][27mm]{0.48\columnwidth}%
\protect\begin{center}
 \begin{tikzpicture}[thick]
 \tikzstyle{every node} = [circle, draw, fill=black!50, 
                           inner sep=1pt, minimum width=3pt]
 \tikzstyle{txt} = [rectangle, draw=none, fill=white, inner sep=-1pt]

 \draw (0,0) to (175:1) node{};
 \draw (0,0) to (185:1) node{};

 \draw (1,0) to ++(-85:1) node{};
 \draw (1,0) to ++(-95:1) node{};

 \draw (1,1) to[bend left]  (2,1) node{};
 \draw (1,1) to[bend right] (2,1) node{};

 \draw (0,0) to ++(-85:1) node{};
 \draw (0,0) to ++(-95:1) node{};

 \draw (0,0) to (1,0) to (1,1) to (0,1) to (0,0);
 \draw (0:0) node{} (1,0) node{} (1,1) node{} (0,1) node{};

 \draw (0.5,0.15) node[txt]{$\lengtha$};
 \draw (0.85,0.5) node[txt]{$\lengthb$};
 \draw (0.5,0.85) node[txt]{$\lengthc$};
 \draw (0.15,0.5) node[txt]{$\lengthd$};

 \draw (-0.5, 0.2) node[txt]{$\lengtha$};
 \draw ( 0.8,-0.5) node[txt]{$\lengthb$};
 \draw ( 1.5, 1.0) node[txt]{$\lengthc$};
 \draw ( 0.2,-0.5) node[txt]{$\lengthd$};

\end{tikzpicture} \begin{tikzpicture}[thick]
 \tikzstyle{every node} = [circle, draw, fill=black!50, 
                           inner sep=1pt, minimum width=3pt]
 \tikzstyle{txt} = [rectangle, draw=none, fill=white, inner sep=-1pt]

 \draw (1,1) to ++( -5:1) node{};
 \draw (1,1) to ++(  5:1) node{};

 \draw (1,0) to ++(-85:1) node{};
 \draw (1,0) to ++(-95:1) node{};

 \draw (0,0) to[bend left]  (-1,0) node{};
 \draw (0,0) to[bend right] (-1,0) node{};

 \draw (0,0) to ++(-85:1) node{};
 \draw (0,0) to ++(-95:1) node{};

 \draw (0,0) to (1,0) to (1,1) to (0,1) to (0,0);
 \draw (0:0) node{} (1,0) node{} (1,1) node{} (0,1) node{};

 \draw (0.5,0.15) node[txt]{$\lengthc$};
 \draw (0.85,0.5) node[txt]{$\lengthb$};
 \draw (0.5,0.85) node[txt]{$\lengtha$};
 \draw (0.15,0.5) node[txt]{$\lengthd$};

 \draw ( 1.5, 0.8) node[txt]{$\lengtha$};
 \draw ( 0.8,-0.5) node[txt]{$\lengthb$};
 \draw (-0.5, 0.0) node[txt]{$\lengthc$};
 \draw ( 0.2,-0.5) node[txt]{$\lengthd$};

\end{tikzpicture}\protect
\par\end{center}%
\end{minipage}\protect
\par\end{centering}

}
\par\end{centering}

\protect\caption{Isospectral quantum graphs with Kirchhoff-Neumann conditions.\label{fig:Isospectral_quantum_graphs}}
\end{figure}
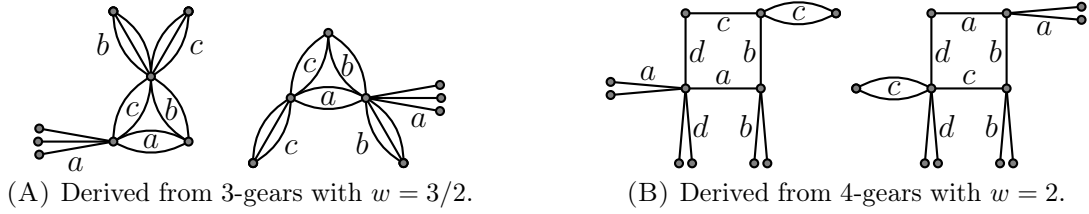

\begin{defn*}
An $\nrTeeth$-gear is a quantum graph with $2\nrTeeth$ vertices
that is comprised of a polygon with $\nrTeeth$ sides of lengths $(\length 1,\length 2,\ldots,\length{\nrTeeth})\in\mathbb{R}_{+}^{\nrTeeth}$
as well as $\nrTeeth$ leaf edges of lengths $(\length 1,\length 2,\ldots,\length{\nrTeeth})$,
called teeth, such that the $i$th tooth is adjacent to the $i$th
polygon side. The dual $\nrTeeth$-gear is obtained by attaching each
tooth at the other vertex of its corresponding polygon side.
\end{defn*}
Figures~\ref{fig:Balanced_three_gears} and~\ref{fig:Unbalanced_three_gears}
show pairs of mutually dual $3$-gears with lengths $(\length 1,\length 2,\length 3)=(\lengtha,\lengthb,\lengthc)$.
We note that an $\nrTeeth$-gear has $\nrTeeth$ leaf vertices and
$\nrTeeth$ polygon vertices, the latter of which have degrees $2$,
$3$, or~$4$. Let $\graph$ be an $\nrTeeth$-gear. We parameterize
the edges of $\graph$ such that its polygon is an oriented cycle
with respect to the induced orientation, and corresponding polygon
sides and teeth are head-to-head or tail-to-tail. Figure~\ref{fig:Polygon_vertices}
indicates the $4$ possible neighborhoods of a polygon vertex. In
order to produce pairs as in Figure~\ref{fig:Isospectral_quantum_graphs},
we introduce an auxiliary weight $\weight>0$ and consider the vertex
conditions described in Figure~\ref{fig:Eigenderivative_transplantation},
where oriented edges and functions on them are denoted by the same
symbol. The corresponding operator is self-adjoint with respect to
the measure obtained by weighting the Lebesgue measures on the teeth
by $\weight$. Namely, if $f$ has restrictions $(\side i)_{i=1}^{\nrTeeth}$
and $(\tooth i)_{i=1}^{\nrTeeth}$, then we consider
\[
\Vert f\Vert_{\weight}^{2}=\langle f,f\rangle_{\weight}=\sum_{i=1}^{\nrTeeth}\int_{0}^{\length i}\side{_{i}}^{2}(x)+\weight\,\tooth i^{2}(x)dx.
\]
If $\varphi$ has restrictions $(\pi_{i})_{i=1}^{n}$ and $(\tau_{i})_{i=1}^{n}$,
then integration by parts gives
\begin{eqnarray*}
\langle\Delta f,\varphi\rangle_{\weight} & = & -\sum_{i=1}^{\nrTeeth}\int_{0}^{\length i}\side{_{i}}''(x)\pi_{i}(x)+\weight\,\tooth i''(x)\tau_{i}(x)dx\\
 & = & \sum_{i=1}^{\nrTeeth}\int_{0}^{\length i}\side{_{i}}'(x)\pi_{i}'(x)+\weight\,\tooth i'(x)\tau_{i}'(x)dx-\sum_{i=1}^{\nrTeeth}\side{_{i}}'(x)\pi_{i}(x)+\weight\,\tooth i'(x)\tau_{i}(x)\Big|_{0}^{\length i}.
\end{eqnarray*}
If $f$ and $\varphi$ obey the vertex conditions described in Figure~\ref{fig:Eigenderivative_transplantation},
then the latter sum vanishes, which can be seen by collecting terms
vertex-by-vertex and considering the $4$ cases shown in Figure~\ref{fig:Polygon_vertices}
separately. In particular, $\langle\Delta f,\varphi\rangle_{\weight}=\langle f,\Delta\varphi\rangle_{\weight}$,
and there exists an $\langle,\rangle_{\weight}$-orthonormal basis
of eigenfunctions with eigenvalue sequence $0=\eigenvalue_{0}<\eigenvalue_{1}\leq\eigenvalue_{2}\leq\ldots$.

We proceed with the core argument, which shows that this sequence
is contained in the spectrum of the corresponding operator on the
dual $\nrTeeth$-gear of $\graph$. If the neighborhood of a polygon
vertex looks like one in the upper half of Figure~\ref{fig:Polygon_vertices},
with edges $\side j$ and $\side k$, then the neighborhood of the
corresponding polygon vertex of the dual $\nrTeeth$-gear looks like
the respective one in the lower half, with edges $\sideOther j$ and
$\sideOther k$, and vice versa. This will allow us to introduce a
linear function, called transplantation, between the spans of non-constant
eigenfunctions on~$\graph$ and those on its dual, which is locally
given by 
\begin{equation}
\sideOther i=\side i'+\weight\,\tooth i'\qquad\text{and}\qquad\toothOther i=\side i'-\tooth i',\qquad\text{or vice versa}.\label{eq:Eigenderivative_transplantation}
\end{equation}
Let $f$ be an eigenfunction on $\graph$ with eigenvalue $\eigenvalue>0$.
Assume that $\graph$ has a polygon vertex of degree $3$ as in the
upper right corner of Figure~\ref{fig:Polygon_vertices}. Let $\side j$,
$\tooth j$, $\side k$, and $\tooth k$ denote the corresponding
restrictions of $f$ so that $\side j''=\eigenvalue\,\side j$, $\tooth j''=\eigenvalue\,\tooth j$,
$\side k''=\eigenvalue\,\side k$, $\tooth k''=\eigenvalue\,\tooth k$,
as well as 
\[
\tooth j'(0)=0,\quad\tooth k'(0)=0,\quad\side j(\length j)=\tooth j(\length j)=\side k(0),\quad\text{and}\quad\side j'(\length j)+\weight\,\tooth j'(\length j)=\side k'(0).
\]
We show that~(\ref{eq:Eigenderivative_transplantation}) gives rise
to a function $\other f$ that is well-defined at the vertex $\sideOther j(\length j)$,
obeys the desired vertex conditions at $\sideOther j(\length j)$
and $\toothOther j(\length j)$, and is an eigenfunction on $\sideOther j$
and $\toothOther j$. Namely,
\begin{equation}
\sideOther j''=\side j'''+\weight\,\tooth j'''=\eigenvalue\,\side j'+\weight\eigenvalue\,\tooth j'=\eigenvalue\,\sideOther j.\label{eq:Transplant_is_eigenfunction}
\end{equation}
Likewise, $\toothOther j''=\eigenvalue\,\toothOther j$, $\sideOther k''=\eigenvalue\,\sideOther k$,
and $\toothOther k''=\eigenvalue\,\toothOther k$. More interestingly,
\begin{eqnarray*}
\toothOther j'(\length j) & = & \side j''(\length j)-\tooth j''(\length j)=\eigenvalue(\side j(\length j)-\tooth j(\length j))=0,\\
\sideOther j(\length j) & = & \side j'(\length j)+\weight\,\tooth j'(\length j)=\side k'(0)=\sideOther k(0)=\toothOther k(0),\\
\sideOther j'(\length j) & = & \eigenvalue(\side j(\length j)+\weight\,\tooth j(\length j))=\eigenvalue(1+\weight)\side k(0)\\
 & = & \eigenvalue(\side k(0)+\weight\,\tooth k(0))+\weight\eigenvalue(\side k(0)-\tooth k(0))=\sideOther k'(0)+\weight\,\toothOther k'(0).
\end{eqnarray*}
Similar arguments apply to the remaining $3$ cases in Figure~\ref{fig:Polygon_vertices}.
In particular, (\ref{eq:Eigenderivative_transplantation})~gives
rise to a globally well-defined eigenfunction $\other f$ on the dual
$\nrTeeth$-gear of $\graph$. For the sake of simplicity, we henceforth
assume that $f$ has restrictions $(\side i)_{i=1}^{\nrTeeth}$ and
$(\tooth i)_{i=1}^{\nrTeeth}$, and that $\other f$ has restrictions
$(\sideOther i)_{i=1}^{\nrTeeth}$ and $(\toothOther i)_{i=1}^{\nrTeeth}$.
The general case just differs by a redistribution of tildes.

It remains to show that for any eigenvalue $\eigenvalue>0$, the transplantation
given by~(\ref{eq:Eigenderivative_transplantation}) restricts to
an injective function between the $\eigenvalue$-eigenspaces. To this
end, we write~(\ref{eq:Eigenderivative_transplantation}) in matrix
form
\[
\left(\begin{array}{c}
\sideOther i\\
\toothOther i
\end{array}\right)=\left(\begin{array}{cc}
1 & \weight\\
1 & -1
\end{array}\right)\left(\begin{array}{c}
\side i'\\
\tooth i'
\end{array}\right),\text{ which leads to }\left(\begin{array}{c}
\side i\\
\tooth i
\end{array}\right)=\frac{1}{\eigenvalue(1+\weight)}\left(\begin{array}{cc}
1 & \weight\\
1 & -1
\end{array}\right)\left(\begin{array}{c}
\sideOther i'\\
\toothOther i'
\end{array}\right).
\]
Since every $\nrTeeth$-gear has $\eigenvalue=0$ as a simple eigenvalue,
we have shown the following.
\begin{thm}
\label{thm:Mutually_dual_gears}Mutually dual $\nrTeeth$-gears are
isospectral for every weight $\weight>0$.
\end{thm}
It is worth mentioning that normalizing~(\ref{eq:Eigenderivative_transplantation})
by the factor $(\eigenvalue(1+\weight))^{-1/2}$ leads to a linear
isometry between the $\eigenvalue$-eigenspaces. More precisely, we
have 
\[
\sideOther i^{2}+\weight\,\toothOther i^{\,2}=(\side i'+\weight\,\tooth i')^{2}+\weight(\side i'-\tooth i')^{2}=(1+\weight)(\side i'^{2}+\weight\,\tooth i'^{2}).
\]
In particular, integration by parts gives
\begin{eqnarray*}
\frac{1}{1+\weight}\Vert\other f\Vert_{\weight}^{2} & = & \sum_{i=1}^{\nrTeeth}\int_{0}^{\length i}\side i'(x)\side i'(x)+\weight\,\tooth i'(x)\tooth i'(x)dx\\
 & = & \sum_{i=1}^{\nrTeeth}\int_{0}^{\length i}\eigenvalue\,\side i^{2}(x)+\weight\,\eigenvalue\,\tooth i^{2}(x)dx+\sum_{i=1}^{\nrTeeth}\side i'(x)\side i(x)+\weight\,\tooth i'(x)\tooth i(x)\Big|_{0}^{\length i}=\eigenvalue\Vert f\Vert_{\weight}^{2}.
\end{eqnarray*}
As above, the latter sum vanishes since $f$ obeys the desired vertex
conditions. Yet, the definition~(\ref{eq:Eigenderivative_transplantation})
has the advantage over its normalized version as it is the same on
all eigenspaces.

Finally, we explain how mutually dual $\nrTeeth$-gears give rise
to pairs as in Figure~\ref{fig:Isospectral_quantum_graphs}. Roughly
speaking, each eigenspace decomposes orthogonally under the action
of the respective graph's isometry group, and the subspace of invariant
elements and its orthogonal complement are transplanted separately.
We exemplify the method with the help of the graphs in Figure~\ref{fig:Derived_graphs_weight_three_halfs}.
Each of them features a $\mathbb{Z}_{2}$-action, given by swapping
parallel edges that make up a side of its central $3$-gon, and a
$\mathbb{Z}_{3}$-action, which moves all but these edges. Thus, each
eigenspace is a unitary representation of $\mathbb{Z}_{2}\times\mathbb{Z}_{3}$.
Eigenfunctions on which $\mathbb{Z}_{3}$ acts by multiplication by
$e^{2\pi i/3}$ or its square vanish on the central $3$-gon. Also,
their derivatives on the $3$ parallel edges of length $\lengthb$
sum to zero at the common polygon vertex, equally for length~$\lengthc$.
For such eigenfunctions, the graphs in Figure~\ref{fig:Derived_graphs_weight_three_halfs}
essentially reduce to the same $3$ subgraphs each of which has $2$
vertices and $3$ edges, and we can tranplant trivially from subgraphs
to subgraphs. Similarly, eigenfunctions that are odd with respect
to the $\mathbb{Z}_{2}$-action are supported on the central $3$-gon,
and we can tranplant trivially from $3$-gon to $3$-gon. Hence, it
suffices to consider the spaces of $\mathbb{Z}_{2}\times\mathbb{Z}_{3}$-invariant
$\eigenvalue$-eigenfunctions on the graphs. However, these spaces
are linearly isometric to the $\eigenvalue$-eigenspaces of the $3$-gears
in Figure~\ref{fig:Unbalanced_three_gears} with $\weight=3/2$,
where the isometry is given by merging the parallel edges of the graphs
in Figure~\ref{fig:Derived_graphs_weight_three_halfs}.

\section{Combinatorial eigenderivative transplantation\label{sec:Combinatorial_eigenderivative_transplantation}}

Long before quantum graphs were introduced as model systems in quantum
chaos~\cite{KottosSmilansky1997,KottosSmilansky1999}, they had been
studied under different names in chemistry and biology, see~\cite{Kuchment2002}
and references therein. Notably, \cite{Below1985} reduces the spectral
analysis of quantum graphs with Kirchhoff-Neumann conditions and commensurable
edge lengths to that of combinatorial graphs with associated row-stochastic
matrices. In contrast to~\cite{KottosSmilansky1999,ShapiraSmilansky2006,BandSmilansky2007,OrenBand2012},
\cite{Below1985} rigorously treats the so-called Dirichlet eigenvalues~\cite{BerkolaikoKuchment2013},
for which (\ref{eq:Eigenfunction_in_terms_of_its_boundary_values})
fails.

In the context of $\nrTeeth$-gears, this leads to a combinatorial
eigenderivative transplantation. We consider pairs of mutually dual
$\nrTeeth$-gears with commensurable edge lengths $(\length 1,\length 2,\ldots,\length{\nrEdges{}})$.
Since scaling all edges by some factor leads to eigenvalues scaled
by that same factor, we may assume that $(\length 1,\length 2,\ldots,\length{\nrEdges{}})\in\mathbb{Z}_{+}^{\nrEdges{}}$.
As is well-known~\cite{BerkolaikoKuchment2013}, adding or removing
a vertex of degree $2$ carrying Kirchhoff-Neumann conditions leaves
the set of eigenfunctions, and therefore the spectrum, unchanged.
We thus replace each edge of length $\length i$ by a path consisting
of $\length i$ edges of length $1$, that is, we subdivide all edges
into edges of unit length.

In the style of~\cite{Below1985}, we first consider Dirichlet eigenvalues
$\eigenvalue=(k\pi)^{2}$ with $k\in\mathbb{Z}$. By~virtue of~(\ref{eq:Boundary_condition_for_Dirichlet_eigenvalues}),
the Neumann conditions at the leaf vertices propagate towards the
polygons so that any $\eigenvalue$-eigenfunction has vanishing first
derivatives not only at the leaf vertices but also at the polygon
vertices in the direction of the respective leaf vertex. In particular,
any $\eigenvalue$-eigenfunction is uniquely determined by its restriction
to the polygon, which is, in fact, a $\eigenvalue$-eigenfunction
of the circle with circumference $\length 1+\length 2+\ldots+\length{\nrEdges{}}$.
On the other hand, any $\eigenvalue$-eigenfunction of this circle
gives rise to unique $\eigenvalue$-eigenfunctions on the $\nrTeeth$-gears.

We turn to eigenvalues $\eigenvalue$ for which~(\ref{eq:Eigenfunction_in_terms_of_its_boundary_values})
holds. In particular, any $\eigenvalue$-eigenfunction $f$ on one
of the subdivided $\nrTeeth$-gears is determined by its values at
the vertices. We differentiate~(\ref{eq:Eigenfunction_in_terms_of_its_boundary_values})
and set $\length i=1$ to obtain
\begin{equation}
\sin(\sqrt{\eigenvalue})f_{i}'(0)=\sqrt{\eigenvalue}\left(-\cos(\sqrt{\eigenvalue})f_{i}(0)+f_{i}(1)\right).\label{eq:Outward_derivative_in_terms_of_values}
\end{equation}
If $\vertex{}$ is a vertex with neighbors $P$ on the polygon and
neighbors $T$ on teeth, then the vertex condition at $\vertex{}$
described in Figure~\ref{fig:Eigenderivative_transplantation} is
given by
\begin{equation}
\cos(\sqrt{\eigenvalue})f(\vertex{})=\frac{1}{|P|+\weight\,|T|}\Bigg(\sum_{\vertex{}'\in P}f(\vertex{}')+\sum_{\vertex{}'\in T}\weight\,f(\vertex{}')\Bigg).\label{eq:Weighted_KN_condition_in_Markov_terms}
\end{equation}
Thus, the values of $f$ at the vertices give rise to a $\cos(\sqrt{\eigenvalue})$-eigenvector
of the row-stochastic matrix $M$ that corresponds to the random walk
on the vertices where edges belonging to teeth are taken $\weight$
times as likely as edges belonging to the polygon. In~fact, the $\eigenvalue$-eigenspace
is isomorphic to the $\cos(\sqrt{\eigenvalue})$-eigenspace of $M$.
For the sake of brevity, we call $M$ the Markov matrix of the $\nrTeeth$-gear,
and denote the right-hand side of~(\ref{eq:Weighted_KN_condition_in_Markov_terms})
by $M[f](\vertex{})$. We note that $M$ is irreducible and has period
$2$ or $1$, depending on whether the subdivided $\nrTeeth$-gear
is bipartite or not. Moreover, $M=D_{1}^{-1}\adjacencyMatrix{}D_{2}$
where $D_{1}$ and $D_{2}$ are invertible diagonal matrices, and
$\adjacencyMatrix{}$ is the adjacency matrix of this graph. In particular,
$M$ is similar to a symmetric matrix and therefore has spectrum in
$[-1,1]$, where $1$ is the simple Perron-Frobenius eigenvalue, and
$-1$ is an eigenvalue precisely if $M$ has period $2$, in which
case it is also simple. In the following, we give an alternative proof
of Theorem~\ref{thm:Mutually_dual_gears}. Since the eigenvalues
of a quantum graph depend continuously on its edge lengths~\cite[Theorem 3.1.2]{BerkolaikoKuchment2013},
it suffices to show the following.
\begin{thm}
The Markov matrices of mutually dual $\nrTeeth$-gears with integral
edge lengths are isospectral for every weight $\weight>0$.\label{thm:Isospectral_Markov_matrices}
\end{thm}
Let $f$ be a function on the vertices of a subdivided $\nrTeeth$-gear.
In view of~(\ref{eq:Outward_derivative_in_terms_of_values}), we
define the outward and inward derivatives of $f$ at the vertex $\vertex{}$
along the edge shared with vertex $\vertex{}'$ as
\[
f'_{[\vertex{},\vertex{}']}(\vertex{})=-M[f](\vertex{})+f(\vertex{}')\qquad\text{and}\qquad f'_{[\vertex{}',\vertex{}]}(\vertex{})=-f'_{[\vertex{},\vertex{}']}(\vertex{})=M[f](\vertex{})-f(\vertex{}').
\]
This definition makes any function satisfy the combinatorial version
of the desired vertex conditions. Namely, if $\vertex{}$ has neighbors
$P$ and $T$ on the polygon and teeth, respectively, then
\begin{equation}
\sum_{\vertex{}'\in P}f'_{[\vertex{},\vertex{}']}(\vertex{})+\sum_{\vertex{}'\in T}\weight\,f'_{[\vertex{},\vertex{}']}(\vertex{})=-(|P|+\weight\,|T|)\,M[f](\vertex{})+\Bigg(\sum_{\vertex{}'\in P}f(\vertex{}')+\sum_{\vertex{}'\in T}\weight\,f(\vertex{}')\Bigg)=0.\label{eq:All_functions_satisfy_combinatorial_vertex_conditions}
\end{equation}
In particular, if $\vertex{}$ has the sole neighbor $\vertex{}'$,
then $f'_{[\vertex{},\vertex{}']}(\vertex{})=0$, and if $\vertex{}$
has degree $2$ and neighbors $\vertex{}'$ and $\vertex{}''$, then
$f'_{[\vertex{}',\vertex{}]}(\vertex{})=f'_{[\vertex{},\vertex{}'']}(\vertex{})$.
For the sake of simplicity, we assume that the underlying $\nrTeeth$-gear
has oriented edges $(\side i)_{i=1}^{\nrTeeth}$ and $(\tooth i)_{i=1}^{\nrTeeth}$
as in Figure~\ref{fig:Polygon_vertices}, the general case is obtained
by redistributing tildes. Similarly to Section~\ref{sec:Eigenderivative_Transplantation},
we denote the restrictions of $f$ to the corresponding previously-introduced
paths by the same symbol, where each of the original $\nrTeeth$ polygon
vertices appears in as many paths as its degree. We orient these paths
as their underlying edges, and define the derivative along $\side i=[\vertex 0,\vertex 1,\ldots,\vertex{\length i}]$
as 
\[
\side i'(\vertex j)=\begin{cases}
f'_{[\vertex j,\vertex{j+1}]}(\vertex j) & \text{if }j<\length i,\\
f'_{[\vertex{j-1},\vertex j]}(\vertex j) & \text{if }j>0,
\end{cases}
\]
similarly for the restrictions $(\tooth i)_{i=1}^{\nrTeeth}$. This
allows to transplant an arbitrary function $f$ via~(\ref{eq:Eigenderivative_transplantation})
to obtain a function $\other f$ with restrictions $(\sideOther i)_{i=1}^{\nrTeeth}$
and $(\toothOther i)_{i=1}^{\nrTeeth}$ on the subdivided dual $\nrTeeth$-gear,
which is well-defined by reason of~(\ref{eq:All_functions_satisfy_combinatorial_vertex_conditions})
and virtually the same argument that showed continuity in the quantum
graph setting. In order to show that this transplantation maps eigenfunctions
of~$M$ to eigenfunctions of its counterpart~$\other M$, assume
that $M[f]=\eigenvalueMarkov f$ for some $\eigenvalueMarkov\in[-1,1]$.
If $[\vertex{}',\vertex{},\vertex{}'']$ is part of one of the paths
$(\side i)_{i=1}^{\nrTeeth}$ or $(\tooth i)_{i=1}^{\nrTeeth}$, say
$\side i$, then
\begin{eqnarray*}
\side i'(\vertex{}')+\side i'(\vertex{}'') & = & -M[f](\vertex{}')+f(\vertex{})-f(\vertex{})+M[f](\vertex{}'')\\
 & = & \eigenvalueMarkov(-f(\vertex{}')+M[f](\vertex{})-M[f](\vertex{})+f(\vertex{}''))=2\eigenvalueMarkov\side i'(\vertex{}).
\end{eqnarray*}
This entails the combinatorial version of~(\ref{eq:Transplant_is_eigenfunction}),
namely, $\other M[\other f](\vertexOther{})=\eigenvalueMarkov\other f(\vertexOther{})$
at all vertices~$\vertexOther{}$ that were introduced with the subdivision.
We therefore turn to vertices of the underlying $\nrTeeth$-gear.
Let $\vertexOther{}$ be the leaf vertex with neighbor $\vertexOther{\toothOther{}}$
on the path $\toothOther j$ in the lower right corner of Figure~\ref{fig:Polygon_vertices},
and $\vertex{}$ be the polygon vertex with neighbors $\vertex{\side{}}$
and $\vertex{\tooth{}}$ on $\side j$ and $\tooth j$ in the upper
one. Then
\begin{eqnarray*}
\other M[\other f](\vertexOther{}) & = & \toothOther j(\vertexOther{\toothOther{}})=\side j'(\vertex{\side{}})-\tooth j'(\vertex{\tooth{}})=-M[f](\vertex{\side{}})+f(\vertex{})-f(\vertex{})+M[f](\vertex{\tooth{}})\\
 & = & \eigenvalueMarkov(-f(\vertex{\side{}})+M[f](\vertex{})-M[f](\vertex{})+f(\vertex{\tooth{}}))=\eigenvalueMarkov(\side j'(\vertex{})-\tooth j'(\vertex{}))=\eigenvalueMarkov\other f(\vertexOther{}).
\end{eqnarray*}
Similarly, let $\vertexOther{}$ and $\vertex{}$ be the polygon vertices
with neighbors $\vertexOther{\tooth k}$, $\vertexOther{\side k}$,
$\vertexOther{\side j}$, $\vertex{\tooth j}$,$\vertex{\side j}$,
and $\vertex{\side k}$ on $\toothOther k$, $\sideOther k$, $\sideOther j$,
$\tooth j$, $\side j$, and $\side k$, respectively, and $\vertex{\tooth k}$
be the neighbor of the leaf vertex on $\tooth k$. Then
\[
\begin{alignedat}{1}(2+\weight)\other M[\other f](\other{\vertex{}}) & =\sideOther j(\vertexOther{\side j})+\sideOther k(\vertexOther{\side k})+\weight\,\toothOther k(\vertexOther{\tooth k})\\
 & =\side j'(\vertex{\side j})+\weight\,\tooth j'(\vertex{\tooth j})+\side k'(\vertex{\side k})+\weight\,\tooth k'(\vertex{\tooth k})+\weight\,\side k'(\vertex{\side k})-\weight\,\tooth k'(\vertex{\tooth k})\\
 & =\eigenvalueMarkov(-f(\vertex{\side j})-\weight f(\vertex{\tooth j})+(1+\weight)f(\vertex{\side k}))+f(\vertex{})+\weight f(\vertex{})-(1+\weight)f(\vertex{})\\
 & =\eigenvalueMarkov(M[f](\vertex{})-f(\vertex{\side j})+\weight(M[f](\vertex{})-f(\vertex{\tooth j}))+(1+\weight)(f(\vertex{\side k})-M[f](\vertex{})))\\
 & =\eigenvalueMarkov(\side j'(\vertex{})+\weight\,\tooth j'(\vertex{})+(1+\weight)\side k'(\vertex{}))\\
 & =\eigenvalueMarkov(\sideOther j(\other{\vertex{}})+\sideOther k(\other{\vertex{}})+\weight\,\toothOther k(\other{\vertex{}}))=(2+\weight)\eigenvalueMarkov\other f(\other{\vertex{}}).
\end{alignedat}
\]
The remaining $3$ cases in Figure~\ref{fig:Polygon_vertices} follow
similarly. Hence, $\other M[\other f]=\eigenvalueMarkov\other f$
whenever $M[f]=\eigenvalueMarkov f$. In order to determine the kernel
of the transplantation, we assume that $\other f=0$. Since $\weight\neq-1$,
we have $\side i'=\tooth i'=0$ for all $i$. Thus, if $[\vertex{},\vertex{}']$
is part of one of the paths $(\side i)_{i=1}^{\nrTeeth}$ or $(\tooth i)_{i=1}^{\nrTeeth}$,
then 
\[
0=f'_{[\vertex{},\vertex{}']}(\vertex{})=-M[f](\vertex{})+f(\vertex{}')=-\eigenvalueMarkov f(\vertex{})+f(\vertex{}'),\text{ which leads to }f(\vertex{}')=\eigenvalueMarkov f(\vertex{}).
\]
Hence, $f(\vertex{})=\eigenvalueMarkov^{\length{}}f(\vertex{})$ for
each vertex $\vertex{}$ on the central polygon of length $\length{}=\length 1+\length 2+\ldots+\length{\nrTeeth}$.
If $f(\vertex{})=0$ for one such vertex $\vertex{}$, then $f$ vanishes
on the entire polygon, and through $M[f]=\eigenvalueMarkov f$ on
the entire subdivided $\nrTeeth$-gear. On the other hand, if $f(\vertex{})\neq0$
for a vertex $\vertex{}$ on the polygon, then $|\eigenvalueMarkov|=1$,
that is, $\eigenvalueMarkov=\pm1$. Hence, the transplantation is
injective on the span of $\eigenvalueMarkov$-eigenvectors with $\eigenvalueMarkov\neq\pm1$.
Since subdivided mutually dual $\nrTeeth$-gears are either both bipartite
or both non-bipartite, we have proven Theorem~\ref{thm:Isospectral_Markov_matrices}.
In order to derive an explicit conjugator for $M$ and $\other M$,
we note that their $\pm1$-eigenspaces are given by functions that
satisfy $f(\vertex{})=\pm f(\vertex{}')$ whenever $\vertex{}$ and
$\vertex{}'$ are neighbors. Thus, if $M[f]=\pm f$, then 
\[
f'_{[\vertex{},\vertex{}']}(\vertex{})=-f'_{[\vertex{}',\vertex{}]}(\vertex{})=-M[f](\vertex{})+f(\vertex{}')=\mp f(\vertex{})+f(\vertex{}')=0.
\]
Hence, the $\pm1$-eigenspaces of $M$ are annihilated by the transplantation.
Let $\allOnesMatrix +$ be a rank-$1$ matrix that maps the $1$-eigenspace
of $M$ onto that of $\other M$, and annihilates all other eigenspaces.
If $-1$ is an eigenvalue of $M$, we choose $\allOnesMatrix -$ analogously,
otherwise we define $\allOnesMatrix -=0$. Writing $M=D_{1}^{-1}\adjacencyMatrix{}D_{2}$
as above, the eigenspaces of $M$ are orthogonal with respect to the
inner product given by the diagonal matrix $D_{1}D_{2}$. In~particular,
we can choose $\allOnesMatrix +=\allOnesMatrix{2\nrTeeth}D_{1}D_{2}$
where $\allOnesMatrix{2\nrTeeth}$ denotes the $2\nrTeeth\times2\nrTeeth$
all-ones matrix. If $T$ denotes the matrix that corresponds to the
transplantation, then $C=T+\allOnesMatrix ++\allOnesMatrix -$ is
invertible and satisfies $\other MC=CM$ by construction. In addition,
the arguments at the end of Section~\ref{sec:Eigenderivative_Transplantation}
equally apply to the subdivided versions of the graphs in Figure~\ref{fig:Isospectral_quantum_graphs},
which shows that their Markov matrices are isospectral.

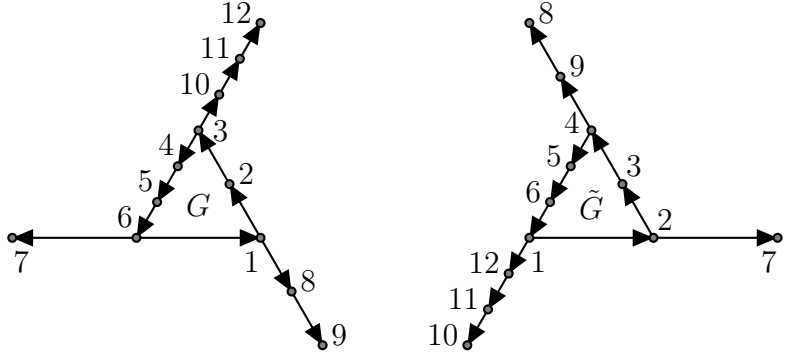
\begin{figure}
\begin{centering}
 \begin{tikzpicture}[scale=0.55, >=triangle 45, thick,
                    decoration={markings,mark=at position 1 with {\arrow{>}}}]
 \tikzstyle{every node} = [circle, draw, fill=black!50, 
                       inner sep=1pt, minimum width=3pt]
 \tikzstyle{txt} = [rectangle, draw=none, fill=white, inner sep=-1pt]

 \draw (1.5,0.8) node[txt]{$\graph$};

 \draw[postaction={decorate}] (0:0) -- ++(180:3);
 \draw[postaction={decorate}] (0:3) -- ++(-60:1.5);
 \draw[postaction={decorate}] (0:3)++(-60:1.5) -- ++(-60:1.5);
 \draw[postaction={decorate}] (60:3) -- ++( 60:1);
 \draw[postaction={decorate}] (60:4) -- ++( 60:1);
 \draw[postaction={decorate}] (60:5) -- ++( 60:1);
 \draw (180:3)          node[txt, below=9pt, right=1pt] {$7$}node{};
 \draw (0:3)++(-60:1.5) node[txt, above=4pt, right=4pt]{$8$} node{};
 \draw (0:3)++(-60:3)   node[txt, above=4pt, right=4pt]{$9$} node{};
 \draw (60:4)           node[txt, above=4pt, left=4pt]{$10$} node{};
 \draw (60:5)           node[txt, above=4pt, left=4pt]{$11$} node{};
 \draw (60:6)           node[txt, above=4pt, left=4pt]{$12$} node{};

 \draw[postaction={decorate}] (0:0) -- (0:3);
 \draw[postaction={decorate}] (0:3) -- ++(120:1.5);
 \draw[postaction={decorate}] (0:3)++(120:1.5) -- (60:3);
 \draw[postaction={decorate}] (60:3) -- (60:2);
 \draw[postaction={decorate}] (60:2) -- (60:1);
 \draw[postaction={decorate}] (60:1) -- (60:0);
 \draw (0:3)            node[txt, below=9pt, left=1pt] {$1$} node{};
 \draw (0:3)++(120:1.5) node[txt, above=4pt, right=4pt]{$2$} node{};
 \draw (60:3)           node[txt, right=6pt]           {$3$} node{};
 \draw (60:2)           node[txt, above=8pt, left=2pt] {$4$} node{};
 \draw (60:1)           node[txt, above=8pt, left=2pt] {$5$} node{};
 \draw ( 0:0)           node[txt, above=8pt, left=2pt] {$6$} node{};
\end{tikzpicture}\hspace*{10mm} \begin{tikzpicture}[scale=0.55, >=triangle 45, thick,
                    decoration={markings,mark=at position 1 with {\arrow{>}}}]
 \tikzstyle{every node} = [circle, draw, fill=black!50, 
                       inner sep=1pt, minimum width=3pt]
 \tikzstyle{txt} = [rectangle, draw=none, fill=white, inner sep=-1pt]

 \draw (1.5,0.8) node[txt]{$\other{\graph}$};

 \draw[postaction={decorate}] (  0:3) -- (0:6);
 \draw[postaction={decorate}] ( 60:3) -- ++(120:1.5);
 \draw[postaction={decorate}] ( 60:3)++(120:1.5) -- ++(120:1.5);
 \draw[postaction={decorate}] (  0:0) -- (240:1);
 \draw[postaction={decorate}] (240:1) -- (240:2);
 \draw[postaction={decorate}] (240:2) -- (240:3);
 \draw (  0:6)            node[txt, below=9pt, left=1pt]{$7$}node{};
 \draw ( 60:3)++(120:1.5) node[txt, above=4pt, right=4pt]{$9$} node{};
 \draw ( 60:3)++(120:3)   node[txt, above=4pt, right=4pt]{$8$} node{};
 \draw (240:1)            node[txt, above=4pt, left=4pt]{$12$} node{};
 \draw (240:2)            node[txt, above=4pt, left=4pt]{$11$} node{};
 \draw (240:3)            node[txt, above=4pt, left=4pt]{$10$} node{};

 \draw[postaction={decorate}] (0:0) -- (0:3);
 \draw[postaction={decorate}] (0:3) -- ++(120:1.5);
 \draw[postaction={decorate}] (0:3)++(120:1.5) -- (60:3);
 \draw[postaction={decorate}] (60:3) -- (60:2);
 \draw[postaction={decorate}] (60:2) -- (60:1);
 \draw[postaction={decorate}] (60:1) -- (60:0);
 \draw (0:3)            node[txt, above=7pt, right=2pt]{$2$} node{};
 \draw (0:3)++(120:1.5) node[txt, above=7pt, right=2pt]{$3$} node{};
 \draw (60:3)           node[txt, above=3pt, left=5pt] {$4$} node{};
 \draw (60:2)           node[txt, above=4pt, left=4pt] {$5$} node{};
 \draw (60:1)           node[txt, above=4pt, left=4pt] {$6$} node{};
 \draw ( 0:0)           node[txt, below=9pt, right=1pt]{$1$} node{};
\end{tikzpicture}
\par\end{centering}

\protect\caption{Zeta-equivalent simple digraphs.\label{fig:Zeta-equivalent-digraphs}}
\end{figure}

Lastly, we mention further presences of the eigenderivative transplantation
method in terms of conjugacy. In fact, the method became apparent
to us when we discovered corresponding conjugators for the matrices
$A^{1}(\sqrt{\eigenvalue})$ and $A^{2}(\sqrt{\eigenvalue})$ in~\cite{OrenBand2012},
which arise when one assumes~(\ref{eq:Eigenfunction_in_terms_of_its_boundary_values})
on all edges, and which characterize eigenvalues through the transcendental
equations $\det(A^{i}(\sqrt{\eigenvalue}))=0$. Another characterization
of the eigenvalues of a quantum graph is given by the scattering approach,
which yields an exact trace formula~\cite{KottosSmilansky1999,BolteEndres2009}.
It can be shown that mutually dual $\nrTeeth$-gears have conjugated
edge $S$-matrices, giving yet another isospectrality proof for the
pairs in Figure~\ref{fig:Isospectral_quantum_graphs}.

Finally, we consider the digraphs $\graph$ and $\other{\graph}$
in Figure~\ref{fig:Zeta-equivalent-digraphs}. Note that the teeth
of $\other G$  are head-to-tail with their corresponding polygon
side. We let $\identityMatrix{12}$ and $\allOnesMatrix{12}$ denote
the $12\times12$ identity and all-ones matrix, respectively. For
$\graph$, we denote its adjacency matrix by $\adjacencyMatrix{\graph}$,
its out-degree matrix by $\outdegreeMatrix{\graph}$, and its in-degree
matrix by $\indegreeMatrix{\graph}$, the latter two of which have
the row sums of $\adjacencyMatrix{\graph}$ and $\adjacencyMatrix{\graph}^{T}$
on their diagonals. For $\boldsymbol{z}=(x,y,\varA,\varAT,\varDout,\varDin)\in\mathbb{C}^{6}$,
we let 
\[
L_{\graph}(\boldsymbol{z})=x\identityMatrix{12}+y\allOnesMatrix{12}+\varA\,\adjacencyMatrix{\graph}+\varAT\,\adjacencyMatrix{\graph}^{T}+\varDout\,\outdegreeMatrix{\graph}+\varDin\,\indegreeMatrix{\graph},
\]
and similarly for $L_{\other{\graph}}(\boldsymbol{z})$. The homogeneous
polynomial $\det(L_{\graph}(\boldsymbol{z}))_{|y=0}$ can be viewed
as a generalized characteristic polynomial. It determines, and is
determined by, the reversing zeta function~\cite{Herbrich2014},
which generalizes the Bartholdi zeta function~\cite{Bartholdi1999}
to digraphs, which in turn generalizes the famous Ihara-Selberg zeta
function~\cite{Ihara1966}. The matrix\setlength{\arraycolsep}{2pt}
\[
T=\left(\begin{array}{cccccccccccc}
\varA^{3} & 0 & 0 & 0 & 0 & 2\varA^{2}\varDout & \varA^{3} & 0 & 0 & 0 & 0 & 0\\
2\varA^{2}\varDout & \varA^{3} & 0 & 0 & 0 & 0 & 0 & \varA^{3} & 0 & 0 & 0 & 0\\
0 & \varA^{2}\varDout & \varA^{3} & 0 & 0 & 0 & 0 & \varA^{2}\varDout & \varA^{3} & 0 & 0 & 0\\
0 & 0 & 2\varA^{2}\varDout & \varA^{3} & 0 & 0 & 0 & 0 & 0 & \varA^{3} & 0 & 0\\
0 & 0 & 0 & \varA^{2}\varDout & \varA^{3} & 0 & 0 & 0 & 0 & \varA^{2}\varDout & \varA^{3} & 0\\
0 & 0 & 0 & 0 & \varA^{2}\varDout & \varA^{3} & 0 & 0 & 0 & 0 & \varA^{2}\varDout & \varA^{3}\\
\varA^{2}\varAT & 0 & 0 & 0 & 0 & 0 & -\varA^{2}\varAT & 0 & 0 & 0 & 0 & 0\\
0 & \varA\varAT^{2} & 0 & 0 & 0 & 0 & 0 & -\varA\varAT^{2} & 0 & 0 & 0 & 0\\
0 & \varA\varAT\varDout & \varA^{2}\varAT & 0 & 0 & 0 & 0 & -\varA\varAT\varDout & -\varA^{2}\varAT & 0 & 0 & 0\\
0 & 0 & 0 & \varAT^{3} & 0 & 0 & 0 & 0 & 0 & -\varAT^{3} & 0 & 0\\
0 & 0 & 0 & \varAT^{2}\varDout & \varA\varAT^{2} & 0 & 0 & 0 & 0 & -\varAT^{2}\varDout & -\varA\varAT^{2} & 0\\
0 & 0 & 0 & 0 & \varA\varAT\varDout & \varA^{2}\varAT & 0 & 0 & 0 & 0 & -\varA\varAT\varDout & -\varA^{2}\varAT
\end{array}\right)
\]
satisfies $L_{\other{\graph}}(\boldsymbol{z})T=TL_{\graph}(\boldsymbol{z})$
and has determinant $((2\varA^{3})^{6}-(2\varA^{2}\varDout)^{6})\varA^{8}\varAT^{10}$,
which can be seen by adding its last $6$ columns to its first $6$
ones. In particular, $\det(L_{\graph}(\boldsymbol{z}))=\det(L_{\other{\graph}}(\boldsymbol{z}))$,
meaning $\graph$ and $\other{\graph}$ are zeta-equivalent and have
zeta-equivalent complements~\cite{Herbrich2014}. The pattern of
non-zero entries in $T$ is reminiscent of combinatorial eigenderivative
transplantation. The conjugator $T$ can be readily generalized to
all pairs coming from dual $\nrTeeth$-gears all of whose polygon
vertices have degree $3$, where teeth have to be oriented as in Figure~\ref{fig:Zeta-equivalent-digraphs}.
However, the graphs in Figure~\ref{fig:Unbalanced_three_gears} do
not lead to zeta-equivalent non-isomorphic digraphs. \setlength{\arraycolsep}{\myArraycolsep}

\bibliographystyle{amsalpha}
\bibliography{changing_gears}

\end{document}